\documentclass[12pt]{article}
\usepackage[russian]{babel}
\usepackage[cp1251]{inputenc}
\usepackage{amsmath,latexsym}
\usepackage[psamsfonts]{amssymb}
\usepackage[dvips]{graphicx}
\usepackage[mathcal]{euscript}
\begin{document}

УДК 517.984
\begin{center}
\textsc{ЭФФЕКТ ЕФИМОВА ОДНОГО МОДЕЛЬНОГО "ТРЕХЧАСТИЧНОГО"
ДИСКРЕТНОГО ОПЕРАТОРА ШРЕДИНГЕРА}
\end{center}
\begin{center}
\textbf{Ю. Х. Эшкабилов}
\end{center}

\begin{abstract}

В работе изучается существование бесконечного числа собственных
значений одного модельного "трехчастичного" оператора Шредингера
Н. Доказывается теорема о необходимом и достаточном условиях для
существования бесконечного числа собственных значений в модели Н,
ниже нижнего края сущестенного спектра.
\end{abstract}

\textit{Ключевые слова:} спектр, существенный спектр,
дискретный \\ спектр, эффект Ефимова.\\

\textbf{1. Введение}\\

Исследование дискретных спектров операторов Шредингера является
наиболее интенсивно изучаемым объектом в теории операторов. Одним
из важных вопросов в спектральном анализе операторов Шредингера
является изучение бесконечности числа собственных значений,
лежащих вне существенного спектра, т.е. существование эффекта
Ефимова в данной модели. Этот эффект впервые был обнаружен
Ефимовым [1] для трехчастичных операторов Шредингера в непрерывном
пространстве. Строгое математическое доказательство существования
эффекта Ефимова было проведено в работе [2], а затем в работах
[3,4] и др.

Помоему, впервые в работе С.Н. Лакаева [5,6] доказано существование
эффекта Ефимова для дискретного оператора Шредингера трехчастичной
системы. Метод был основан на аналитических свойствах определителей
Фредгольма. Затем в работах [7-10] изучено существование эффекта
Ефимова для трехчастичных дискретных операторов Шредингера методом
А.В. Соболева (т.е. методом асимптотики) [11].

В настоящей работе изучено существование эффекта Ефимова для одного
модельного дискретного "трехчастичного" оператора Шредингера,
возникающих в модели Хаббарда. При исследовании эффекта Ефимова, мы
будем пользоваться инструментами принципа минимакса для ограниченных
самосопряженных операторов и свойствами положительных интегральных
операторов. В статье, главным образом приведены достаточные и
необходимые условия для существавания эффекта Ефимова, ниже нижнего
края существенного спектра в данной модели. \\

\textbf{2. Некоторые обозначения и необходимые сведения}\\

Пусть  $\Omega_1\subset \mathbb R^{d_1}, d_1\in N $ и
$\Omega_2\subset\mathbb R^{d_2}, d_2\in N $ - компактные множества.
Рассмотрим непрерывные функции $k_1(x,s)$ на $\Omega^2_1,\,\
k_1(x,s)=\overline{k_1(s,x)}$ и $k_2(y,t)$ на $\Omega^2_2,\,\
k_2(y,t)=\overline{k_2(t,y)}$. В гильбертовом пространстве
$L_2(\Omega_1 \times \Omega_2)$ рассмотрим частично интегральные
операторы $T_1$ и $T_2:$

$$T_1 f(x,y)=\int\limits_{\Omega_1}{k_1(x,s)f(s,y)d\nu_1(s)},\,\
 f\in L_2(\Omega_1 \times \Omega_2),$$
$$T_2 f(x,y)=\int\limits_{\Omega_2}{k_2(y,t)f(x,t)d\nu_2(t)},\,\
f\in L_2(\Omega_1 \times \Omega_2).$$

Пусть $k_0(s,x)$ - произвольная вещественнозначная непрерывная
функция на $\Omega_1 \times \Omega_2.$ Обозначим через $H_0$
оператор умножения на функцию $k_0(x,y),$ т.е. $$H_0 f(x,y)=k_0
(x,y) f(x,y),\,\,\ f \in L_2 (\Omega_1 \times \Omega_2).$$

Рассмотрим линейный ограниченный самосопряженный оператор
$$H=H_0-T \eqno (1)$$ действующий в пространстве $L_2
(\Omega_1 \times \Omega_2),$ где $$T=T_1+T_2. \eqno (2)$$

Оператор $H \,\ (1)$ является общим моделом "трёхчастичного"
дискретного оператора Шредингера, возникающий в модели Хаббарда на
примесной решетке [см. 12].

Сушественный спектр оператора $H$ в более общем виде изучен в
работе [12]. В работе [13] изучен существенный и дискретный спектр
операторов в модели (1), когда потенциальная функция $k_0 (x,y)$
имеет специальный вид: $k_0 (x,y)=u(x)+h(y).$

Через $\rho(\cdot),\,\ \sigma(\cdot),\,\ \sigma_e (\cdot)$ и
$\sigma_d(\cdot)$ обозначим, соответственно, резольвентное
множество, спектр, существенный спектр и дискретный спектр
самосопряженных операторов [14].

В монографии [15] с помощью принципа минимакса изучены
спектральные свойства заданных самосопряженных операторов
ограниченных снизу, действующих в сепарабельном гильбертовом
пространстве. В частности, доказано существование собственных
значений и бесконечность дискретного спектра (т.е. существование
эффекта Ефимова) для некоторых многочастичных гамильтонианов.
Сначала, мы изложим принцип "минимакса" для ограниченных
самосопряженных операторов, которым будем пользоватся при
исследовании на конечность или бесконечность дискретного спектра в
модели (1).

Пусть $ \mathcal{H}-$ сепарабельное гильбертово пространство,
$A:\mathcal{H}\rightarrow \mathcal{H}-$ линейный ограниченный
самосопряженный оператор.

Положим $$\emph{E}_{min}=\emph{E}_{min}(A)=inf \{ \lambda :
\lambda \in \sigma _e (A)\}.$$ Имеем $\emph{E}_{min}(A) \in \sigma
_e (A).$ Число $\emph{E}_{min}(A)$ будем называть нижным краем
существенного спектра оператора $A.$

Определим вещественное число
$$\emph{S}_{min}=\emph{S}_{min}(A)=inf \{(Ax,x): \,\ x \in
\mathcal{H},\,\ \|x\|=1 \}.$$ Тогда $\emph{S}_{min}(A) \in \sigma
 (A)$ и существует элемент $x_1 \in \mathcal{H},\,\ \|x_1\|=1$
такой, что $\emph{S}_{min}(A)=(Ax_1, x_1).$

 Построим ограниченную возрастающую последовательность
$\mu_n=\mu_n(A), \\ n \in N$ следующим способом: $$\mu_1 (A) =
\emph{S}_{min}(A)=(A x_1,x_1)$$ и
$$\mu^{(1)}_{\kappa}(A)=\mu_1(A_{\kappa-1})=\emph{S}_{min}
(A_{\kappa-1}),\,\ \kappa=2,3,...$$ где $A_{\kappa}-$ сужение
оператора $A$ на подпространство
$$\mathcal{H}_{\kappa}=\{x \in \mathcal{H}: \,\ (x,x_i)=0 ,\,\
i=1,2,...,\kappa\},$$ Здесь
$\emph{S}_{min}(A_{\kappa-1})=(A_{\kappa}x_{\kappa}, x_{\kappa}),\,\
x_{\kappa} \in \mathcal{H}_{\kappa-1},\,\ \|x_{\kappa}\|=1,\,\
\kappa=2,3,... \,\ .$

Положим  $$\mu_n (A)= sup \{\mu_1 (A),\,\mu^{(1)}_2(A),\,\ ...,
\mu^{(1)}_n (A)\},\,\ n\in N. \eqno (3)$$ \\
\textbf{Теорема 2.1.} \emph{(принцип минимакса для ограниченных
операторов в операторной форме). Пусть $A:\mathcal{H}\rightarrow
\mathcal{H}-$ ограниченный самосопряжённый оператор. Тогда для
каждого фиксированного $n \in N,$ \\
либо \\
(а) существует $n$ собственных значений (считая вырожденные
собственные значения столько раз,какова их кратность), лежащих ниже
нижнего края $\emph{E}_{min}(A)$ существенного спектра, а
$\mu_{n}(A)$ \emph{из (3) есть} $n-$ е собственное значение
оператора $A$ (с учетом кратности),\\
либо\\
(б) $\mu_n (A)-$ нижний край существенного спектра оператора $A,$
т.е. $\mu_n (A)= \emph{E}_{min}(A),$ при этом $\mu_n (A)=
\mu_{n+\kappa}(A), \forall \kappa \in N$ и существует самое большее
$n-1$ собственных значений (с учетом кратности), лежащих ниже
нижнего края $\emph{E}_{min}(A)$ существенного спектра оператора}
$A.$ \\

Линейный ограниченный оператор $A: \mathcal{H}\rightarrow
\mathcal{H}$ называется положительным и пишется $A\geq 0 $ или
$0\leq A,$ если $(Ax, x )\geq 0,\,\ \forall x \in
\mathcal{H}.$  \\

\textbf{Cледствие 2.1.} \emph{Пусть $A,B: \mathcal{H}\rightarrow
\mathcal{H} -$ линейные ограниченные самосопряженные операторы,
$A\leq B \, (B-A\geq 0)$ и $\emph{E}_{min} (A)=\emph{E}_{min}
(B).$ Тогда }
$$\mu_n(A)\leq \mu_n (B),\ n \in N. \eqno (4)$$ \\

\textbf{3. Вспомогательные утверждения} \\

Изучим конечность и бесконечность дискретного и существенного
спектра самосопряженного оператора $T$ (2). Определим компактные
интегральные операторы  $K_1$ и $K_2$,  действующие соответственно
в $L_2 (\Omega_1)$ и $L_2 (\Omega_2)$ по правилам
$$ (K_1 \varphi)(x)=\int \limits_{\Omega_1}{ k_{1}(x,s)
\varphi(s)d \nu_{1}(s)},$$
$$(K_{2}\phi)(y)=\int \limits_{\Omega_{2}}{ k_{2}(y,t) \phi(t)d
\nu_{2}(t)}.$$ Оператор $T$ является унитарно эквивалентным
оператору $K=K_1\bigotimes E+E \bigotimes K_2$ действующего в
пространстве $L_2 (\Omega_1)\bigotimes L_2
(\Omega_2)$(см.[12,16]), где $E$ - тождественный оператор,
$\bigotimes$ - тензорное произведение. Следовательно, получим [14]
$\sigma (T)=\sigma (K_1)+\sigma (K_2).$ Отсюда имеем
[16]$$\sigma_e (T)=\sigma(K_1)\cup\sigma(K_2)=\{0\}\cup
\sigma_d(K_1)\cup \sigma_d(K_2),$$ $$\sigma_d(T)=\{\omega :
\omega=\alpha+\beta\notin\sigma(K_1)\cup \sigma(K_2),\ \alpha \in
\sigma_d(K_1),\ \beta \in \sigma_d(K_2)\}.$$

Легко заметить, что мощность множеств $\sigma_e (T)$ и $\sigma_d
(T)$ не более чем счетна. Из представления существенного и
дискретного спектра оператора $T,$ легко следует \\

\textbf{Предложение 3.1.} \emph{Имеют место следующие неравенства: \\
(а)$|\sigma_d (K_{\kappa})|+1\leq|\sigma_e(T)|,\ \kappa=1,2;$ \\
(б)$|\sigma_e (T)|\leq|\sigma_d(K_1)|+|\sigma_d(K_2)|+1;$ \\ (в)
$|\sigma_d (T)|\leq|\sigma_d(K_1)|\cdot |\sigma_d(K_2)|,$ \\
где $|\cdot|-$ мощность множества.} \\

\textbf{Предложение 3.2.}\emph{ Дискретный спектр (существенный
спектр) оператора $T$ конечен, тогда и только тогда, когда
дискретный спектр компактных операторов $K_1$ и $K_2$ конечен.}

\emph{Доказательство.} Достаточность очевидна. Докажем
необходимость. Допустим обратное: дискретный спектр оператора $T$
конечен, а дискретный спектр оператора $K_1$ или $K_2$ бесконечен.
Надо подчеркнуть, что дискретный спектр самосопряженных компактных
операторов, состоит из всех их собственных значений отличных от
нуля. Не нарущая общности мы считаем, что $\sigma_d(K_1)-$ конечен
и $\sigma_d(K_2)-$ бесконечен. Тогда существует последовательность
$\beta_n,\ n \in N$ собственных значений оператора $K_2$ такая,
что $\lim\limits_{n\to\infty} \beta_n =0.$

Пусть $\alpha\in \sigma_d (K_1).$ Тогда $\alpha+\beta_n \in
\sigma(T),\,\ n \in N.$ Будем выбрасывать из членов
последовательность $\omega_n= \alpha+ \beta_n,\ n \in N$ те члены,
для которых $\omega_\kappa\in \sigma_d(K_1).$ Тогда количество
выбрашенных членов $ \omega^{\prime}_1,\ \omega^{\prime}_2,\,\ ...
,\ $ из $ \{\omega_n\}_{n \in N}$ конечно, так как само множество
$\sigma_d(K_1)$ конечно. Обозначим через
$\{{\omega}^{(1)}_{\kappa}\}_{\kappa \in N}$ подпоследовательность
из $\{\omega_n\}_{ n \in N}$ таких, что
${\omega}^{(1)}_{\kappa}\notin\sigma_d(K_1),\ \forall \kappa \in
N.$ Теперь будем выбрасывать из
$\{{\omega}^{(1)}_{\kappa}\}_{\kappa \in N}$ те члены, для которых
${\omega}^{(1)}_{\kappa}\in\sigma_d(K_2).$ Однако, количество
выбрашенных членов $\omega^{\prime\prime}_1,\
\omega^{\prime\prime}_2,\  ... ,\ $ из
$\{{\omega}^{(1)}_{\kappa}\}_{\kappa \in N}$ будет конечно, иначе
существовала бы последовательность
$\{\omega^{\prime\prime}_{\kappa}\}_{\kappa\in N}$ собственных
значений компактного оператора $K_2$ такая, что
$\lim\limits_{\kappa\to\infty} \omega^{\prime\prime}_{\kappa}
=\alpha\ne0,$ а это противоречит компактности оператора $K_2.$

Обозначим через $\{{\omega}^{(2)}_{\kappa}\}_{\kappa \in N}$
последовательность состоящую из членов последовательность
$\{{\omega}^{(1)}_{\kappa}\}_{\kappa \in N}$ такую, что
${\omega}^{(2)}_{\kappa} \notin \sigma_d(K_2),\,\ \forall \kappa
\in N.$ Следовательно, имеем
${\omega}^{(2)}_k=\alpha+\beta_{n_{k}}\in \sigma_d(T),\,\ \forall
\kappa \in N,$ где $\beta_{n_{\kappa}}-$ подпоследовательность
последовательности $\{\beta_n\}_{n \in N}\subset \sigma_{d}(K_2).$
Отсюда, $|\sigma_d(T)|=\infty,$ а это противоречит
предположению.$\blacksquare$ \\

\textbf{Предложение 3.3.} \emph{Дискретный спектр (существенный
спектр) оператора $T$ бесконечен, тогда и только тогда, когда
дискретный спектр оператора $K_1$ или $K_2$ бесконечен.}
 \\

\textbf{4. Эффект Ефимова} \\

Существенный спектр $\sigma_e(H)$ оператора $H$ представляется в
виде (Теорема 3.3 [12]):
$$\sigma_e(H)=\sigma(H_0)\cup \sigma(H_0-T_1)
\cup\sigma(H_0-T_2).$$

В дальнейшем, предположим что в модели (1) $\,\ k_0(x,y)-$
произвольная неотрицательная функция на $\Omega_1\times \Omega_2$ и
$k_0^{-1}(\{0\})\cap\Omega_1\times \Omega_2 \neq {\O}$ и $ K_1\geq
0,\,\ K_2\geq 0.$

\textbf{Теорема 4.1.} \emph{Пусть
$H_0\geq(\emph{E}_{min}(H)+\eta_0) E,$  где $\eta_0=sup \sigma_e
(T).$ Если дискретный спектр $\sigma_d (T)$ оператора $T\,\ (2)$
конечен, тогда количество собственных значений в модели (1),
лежащих ниже нижнего края} $\emph{E}_{min}(H)$ \emph{сушественного
спектра $\sigma_e (H)$
конечен.} \\

\emph{Доказательство.}  Имеем
$$\sigma_e (\eta_0 E-T)=\{\omega: \omega=\eta_0-\lambda,\,\ \lambda \in \sigma_e (T)\}.$$
Тогда  $\sigma_e (\eta_0 E-T)\subset [0,\infty)$ и  $0 \in
\sigma_e (\eta_0 E-T ).$ Отсюда
$\emph{E}_{min}(T)=\emph{E}_{min}(\eta_0 E-T)=0.$ Тогда получим,
что
$$\emph{E}_{min}(H-\emph{E}_{min}(H))=\emph{E}_{min}(\eta_0 E-T)=0.$$
Из $H_0\geq(\emph{E}_{min}(H)+\eta_0) E$ следует, что $\eta_0
E-T\leq H-\emph{E}_{min}(H)\cdot E.$ Тогда в силу неравенства (4)
$$\mu_{\kappa}(\eta_0 E-T)\leq\mu_{\kappa}(H-\emph{E}_{min}(H)\cdot E),\,\ \kappa\in N.$$

Пусть $|\sigma_d(T)|=m.$ Тогда для множества
$$\sigma_{d}^{0}(\eta_0 E-T)=\{\omega: \omega\in \sigma_d
(\eta_0 E-T),\,\ \omega<0\}$$ имеем $$1\leq|\sigma_{d}^{0}(\eta_0
E-T)|=m_0\leq m,$$ т.е. число собственных значений оператора $\eta_0
E-T,$ лежащих ниже нижнего края существенного спектра
$\sigma_{e}(\eta_0 E-T)$ не больше чем $m.$

Следовательно, в силу теоремы 2.1 существует натуральное число $n$
такое, что
$$\mu_{n+\kappa}(\eta_0 E-T)=0,\ \forall \kappa \in N \cup \{0\}.$$
Отсюда, имеем $\mu_{n+\kappa}(H-\emph{E}_{min}(H)\cdot E)=0,\,\
\forall \kappa\in N \cup \{0\},$ т.е. количество отрицательных
собственных значений оператора $H-\emph{E}_{min}(H)\cdot E$
конечно. Следовательно, количество собственных значений оператора
$H$ лежащих ниже нижнего края $\emph{E}_{min}(H)$ существенного
спектра, также будет конечной.$\blacksquare$ \\

Из теоремы 4.1. следуют \\

\textbf{Теорема 4.2} \emph{Пусть
$H_0\geq(\emph{E}_{min}(H)+\eta_0) E.$  Чтобы в модели (1) ниже
нижнего края $\emph{E}_{min}(H)$ существенного спектра
$\sigma_e(H)$ существовал эффект Ефимова необходимо, чтобы
оператор $T$ имел бесконечный дискретный спектр, т.е.}
$|\sigma_d(T)|=\infty.$ \\

В пространстве $L_2(\Omega_1 \times \Omega_2)$  определим
самосопряженный операторы $W_1$ и $W_2:$

$$W_1=H_0-T_1,\,\ W_2=H_0-T_2.$$ \\
Надо отметить,  что дискретный спектр у оператора $W_1$ и $W_2$
отсутствует. Так как, если $f (x,y)$ собственная функция оператора
$W_1,$ т.е. $W_1 f=\lambda f,$ для некоторого $\lambda \in
\mathbb{R},$ то для функций $g (x,y)=\psi (y) f (x,y)$ имеем $W_1
g (x,y)=\lambda g (x,y),$ где $\psi (y)$ - произвольная
существенно ограниченная измеримая функция заданная на $\Omega_2.$
Следовательно, всякое собственное значение оператора $W_1,$
является бесконечнократным, т.е. $\sigma_d(W_1)=\emptyset.$
Аналогично можно показать, что $\sigma_d(W_2)=\emptyset.$

    В силу теоремы 3.3 из [12] вытекает, что для нижнего края $\emph{E}_{min}(H)$
   существенного спектра оператора $H$ имеется два случае: $\emph{E}_{min}(H)=\emph{E}_{min}(W_1)$ или
$\emph{E}_{min}(H)=\emph{E}_{min}(W_2).$
\\

 \textbf{Теорема 4.3.} \emph{Пусть
$\emph{E}_{min}(H)=\emph{E}_{min}(W_1)=\Lambda_1\leq
\emph{E}_{min}(W_2).$ Если существует ортонормированная система
$\{\varphi_n\}_{n \in N}\subset L_2(\Omega_2)$ удовлетворяющая
условию
$$(W_1 f_{\kappa}, f_{\kappa})<\Lambda_1+(T_2 f_{\kappa}, f_{\kappa}),\ \kappa\in N,
\eqno (5)$$ где $f_{\kappa}(x,y)=\varphi_0(x)
\varphi_{\kappa}(y),\varphi_0(x)\equiv1,$
 тогда ниже нижнего края существенного
спектра оператора $H$ существует эффект Ефимова.}

\emph{Доказательство.}  Пусть в пространстве $L_2 (\Omega_2)$
существует ортонормированная система $\{\varphi_n\}_{n \in N}$
удовлетворяющая условию (5). Легко проверить, что $$\|T_2
f_n\|=\sqrt{\nu_1 (\Omega_1)}
 \|K_2 \varphi_n\|,\ n \in N.$$ В силу компактности оператора $K_2$ имеем
       $\lim_{n \rightarrow \infty} \|K_2 \varphi_n\|=0.$ Однако, имеем
        $$\sigma(W_1)=\sigma_e(W_1).$$ Отсюда и  из неравенство (5)
получим, что

$$\Lambda_1\leq (W_1 f_n, f_n)<\Lambda_1+(T_2 f_n, f_n),\ n\in N.$$
Следовательно  имеем  $\lim_{n  \rightarrow \infty} (W_1 f_n,
f_n)=\Lambda_1.$ Значит  $\lim_{n  \rightarrow \infty} (H f_n,
f_n)=\Lambda_1.$ С другой стороны $(H f_n, f_n)<\Lambda_1,\ n\in
N.$ Тогда, не нарущая общности, мы  можем предпологать, что
система $\{f_n\}_{n \in N}$ "упорядочена" в следующем смысле:
$$(H f_{\kappa},\ f_{\kappa})\leq (H f_{\kappa+1}, f_{\kappa +
1}),\ \kappa\in N.$$ Имеем $$\emph{S}_{min}(H)=\mu_1(H)\leq (H
f_1, f_1)<\Lambda_1,$$ т.е. сушествует $g_1 \in
L_2(\Omega_1\times\Omega_2), \|g_1\|=1$ такая, что $\mu_1(H)=(H
g_1, g_1)<\Lambda_1$ и по теореме о принципе минимакса число
$\mu_1(H)$ является собственным значением оператора $H.$

Для каждого $\kappa \in N$ определим подпространство
$L_{\kappa}\subset L_2(\Omega_1\times\Omega_2):$
$$L_{\kappa}=\{f\in L_2(\Omega_1\times\Omega_2) : (f, f_j)=0,\
j=1, 2, ..., \kappa \}.$$

Пусть $H_{\kappa}$ сужение оператора $H$ на подпрастранство
$L_{\kappa}, \kappa \in N.$ Для оператора $H_1$ имеем
$$\mu_1(H_1)=\emph{S}_{min}(H_1)=inf\{(H_1 f, f): f \in L_1,
\|f\|=1 \} \leq (H_1 f_2, f_2)=(H f_2, f_2)<\Lambda_1.$$

Так как, $f_2 \in L_1, \|f\|=1,$ то существует $g_2 \in L_1,
\|g_2\|=1 $ такая, что $\mu_1(H_1)=(H_1 g_2, g_2)=(H g_2,
g_2)<\Lambda_1.$

Для каждого $\kappa\geq 2$ повторяя анологичное рассуждение, имеем
$$\mu_1(H_{\kappa})=\emph{S}_{min}(H_{\kappa})=inf\{(H_{\kappa} f,
f): f \in L_{\kappa}, \|f\|=1\}\leq $$ $$ \leq(H_{\kappa}
f_{\kappa+1}, f_{\kappa+1})=(H f_{\kappa+1}, f
_{\kappa+1})<\Lambda_1.$$ Посколько, $f_{\kappa+1}\in L_{\kappa},
\|f_{\kappa+1}\|=1,\,\ \kappa\geq 2.$ Тогда существует
$g_{\kappa+1}\in L_{\kappa},\,\ \|g_{\kappa+1}\|=1$ такая, что
$$\mu_1(H_{\kappa})=(H_{\kappa}g_{\kappa+1}, g_{\kappa+1})=(H
g_{\kappa+1}, g_{\kappa+1})<\Lambda_1.$$

Таким образом, для каждого $\kappa\in N $ оператор сужения
$H_{\kappa}$ имеет собственное значение
$$\mu_1(H_{\kappa})<\Lambda_1\,\ \mbox{и} \,\,\,
\mu_1(H_{\kappa})\leq\mu_1(H_{\kappa+1}),\ \kappa\in N. $$

Следовательно, каждое число
$\omega_{\kappa}=\mu_1(H_{\kappa})<\Lambda_1,\ \kappa=0, 1, 2,
...,$ является собственным значением оператора $H,$ где $H_0=H$ и
$\lim\limits_{n\to\infty} { \omega_n}=\Lambda_1.$
$\blacksquare$ \\

\textbf{Теорема 4.4.} \emph{Пусть
$\emph{E}_{min}(H)=\emph{E}_{min}(W_2)=\Lambda_2\leq
\emph{E}_{min}(W_1).$ Если существует ортонормированная система
$\{\varphi_n\}_{n \in N}\subset L_2(\Omega_1)$ удовлетворяющая
условию
$$(W_2 f_{\kappa}, f_{\kappa})<\Lambda_2+(T_1 f_{\kappa}, f_{\kappa}),\ \kappa\in N,
\eqno (6)$$ где $f_{\kappa}(x,y)=\varphi_{\kappa}(x)\varphi_0(y)
,\varphi_0(y)\equiv1,$
 тогда ниже нижнего края существенного
спектра оператора $H$ существует эффект Ефимова.}\\

\textbf{5. Пример} \\

Рассмотрим последовательность  $p_0=0,\,\ p_1=\frac{{1}}{2},\
p_n=p_{n-1}+\frac{{1}}{2^n},\ n \in N.$ Положим
$$q_n=\frac{{p_n-p_{n-1}}}{2},\ n\in N. $$\\
 На $[0,1]$  определим
функцию $u(x)$
$$ u(x)= \left\{\begin{array}{cc}
  0,\ & \mbox{если} \, x \in [0, \frac{{1}}{2}]  \\
   u_0(x),\ & \mbox{если} \,  x \notin[0, \frac{{1}}{2}]  \\
\end{array}\right.
$$
где $u_0(x)= \sum\limits_{n\in N}{\delta_n r_n (x),}$
$$r_\kappa (x)= \left\{ \begin{array}{ll}
  \frac {{p_{\kappa-1}-x }}{p_{\kappa-1}-q_\kappa},
  \ & \mbox{если} \, x \in [p_{\kappa-1}, q_\kappa]  \\
\frac{{p_\kappa-x}}{p_\kappa-q_\kappa},
\ & \mbox{если} \, x \in [q_\kappa,p_\kappa]  \\
 \end{array}\right.$$
  и $$\delta_1= 1, \,\ \delta_n \leq( \frac{{\sqrt{2}}}{3})^{n},\ n\geq2.$$
   В пространстве $L_2[0,1]$ рассмотрим последовательность ортонормированных функций
$$\varphi_{n}(y)=2^{\frac{{n+1}}{2}}\sin\xi_{n}(y),$$ где

$$ \xi_{\kappa}(y)= \left\{\begin{array}{cc}
  \frac{{\pi}}{p_\kappa-p_{\kappa-1}}
(y-p_{\kappa-1}),\ & \mbox{если} \, y \in [p_{\kappa-1}, p_\kappa]  \\
   0,\ & \mbox{если} \,  y \notin[p_{\kappa-1}, p_\kappa] .  \\
\end{array}\right.
$$

Определим ядро $k_2(y,t),$ заданной  формулой $$k_2(y,t)=
\sum\limits_{n\in N}{(\frac{{2}}{3})^{n} \varphi_n (y)
\varphi_n(t).} \eqno (7)$$ Ряд (7) равномерно сходится на квадрате
$[0, 1]^2.$ Следовательно, интегральный оператор $K_2,$ заданный
ядром $k_2(y,t),$ является самосопряженным и положительным в
$L_2[0, 1].$

В пространстве $L_2{[0,1]}$ рассмотрим следующую модель
$$H=H_0-(\gamma T_1+T_2),\,\  \gamma\geq\frac{{2}}{3} \eqno (8)$$ где $$H_0 f(x,y)=u(x) u(y) f(x,y)$$
$$T_1 f(x,y)= \int\limits_0^1{f(s,y)d\nu(s)},\,\ $$
$$T_2 f(x,y)=\int\limits_0^1{k_2(y,t) f(x,t)d\nu(t)}.$$

Покажем, что в модели (8) ниже нижнего края существенного спектра
существует эффект Ефимова. Раcсмотрим последовательность
ортонормированных функций $f_n\in L_2 [0,1]^2,\,\ n\in N:$
$$f_{n}(x,y)=\varphi_0(x) \varphi_{n} (y),\,\ n \in N,$$ где
$\varphi_0(x)\equiv1.$ Имеем

$$((H_0-\gamma T_1)f_1, f_1)=(H_0 f_1, f_1)-\gamma(T_1 f_1, f_1)=-\gamma(T_1 f_1, f_1)=-\gamma.$$
Отсюда и из неравенство

$$((H_0-\gamma T_1)f, f)\geq-\gamma(T_1 f, f)\geq-\gamma,\,\ f\in L_2 [0,1]^2$$
получим, что $\emph{E}_{min}(H_0-\gamma T_1)=-\gamma.$ С другой
стороны, имеем
$$((H_0-T_2)f, f)=(H_0 f,f)-(T_2 f,f)\geq -(T_2 f, f)\geq -{\frac{{2}}{3}},\,\ f\in L_2 [0,1]^2 .$$
Следовательно, из $\gamma\geq\frac{{2}}{3}$  и в силу теоремы 3.3
из [12] получим, что $\emph{E}_{min}(H_0-(\gamma T_1+
T_2))=-\gamma.$ Положим $b_n=(\frac{{2}}{3})^{n},\,\ n \in N.$
Тогда каждое число $\omega_n=\gamma+b_n$ является собственным
значением оператора $T=\gamma T_1+T_2$ так, как $$(T f_n)(x,
y)=\gamma\int\limits_0^1{\varphi_0(s)
\varphi_n(y)d\nu(s)}+\int\limits_0^1{k_2(y,t)\varphi_0(x)\varphi_n(t)d\nu(t)}=\gamma
f_n(x, y)+$$ $$+b_n f_n(x, y)= (\gamma+b_n)f_n(x,y),\,\ n \in N.$$ \\

Для оператора $H_0$ получим, что
$$(H_0 f_n, f_n)=\int\limits_0^1{u (s)
d\nu(s)}\int\limits_0^1{u (t)\varphi_n^2
(t)d\nu(t)}\leq\int\limits_0^1{u (t)\varphi_n^2
(t)d\nu(t)}\leq\delta_n,\,\ n\in N.$$ Отсюда,

$$((H_0-\gamma T_1)f_n, f_n)\leq\delta_n -\gamma<-\gamma+(\frac{{2}}{3})^{n}=
\emph{E}_{min}(H)+(T_2 f_n, f_n),\,\ n\geq2 $$ т.е. выполняется
условие из теоремы 4.3.  Значит, в модели (8) ниже нижнего края
сущесвенного спектра существует эффект Ефимова.

\newpage

\begin{center} \textbf{Сведения о авторе} \end{center}

Ф.И.О: Эшкабилов Юсуп Халбаевич.

Место работы: механико-математический факультет Национального
Университета Узбекистана им М.Улугбека.

Адрес по месту работы: 100174, Узбекистан, г Ташкент, ВУЗгородок.

Домашний адрес: 100139, Узбекистан, г Ташкент, Чиланзар 23-34-42.

Тел: 8-371-246-02-30 (раб);
     8-371-274-67-10 (дом);

e-mail: yusup62@rambler.ru

\end{document}